\documentclass[12pt,a4]{article}
\usepackage{graphicx,amssymb,amsmath,amsthm,cite}
\newtheorem{theorem}{Theorem}

\newtheorem{remark}{Remark}

\def\Z{\mathbb{Z}} 
 
\def\D{{\cal D}}%^{[c,d]}}
\def\Dch{{\cal D}}%^{[0,8]}}
\def\L{{L^2[c,d]}} 
\def\V{V}%^{[c,d]}}
\def\W{W}%^{[c,d]}}
\def\dicv{{\cal V}}%^{[c,d]}}
\def\dicw{{\cal W}}%^{[c,d]}}
\def\kj{k_\ell} 
\def\kl{p(n)} 
\def\kr{p(n)} 
\def\la{\langle}
\def\ra{\rangle}

\DeclareMathOperator{\supp}{supp} 
\DeclareMathOperator{\Span}{span}

\title{From cardinal spline wavelet bases  
to highly coherent dictionaries}
\author{Miroslav ANDRLE and Laura REBOLLO-NEIRA\\%\thanks{Support from EPSRC (GR$/$R86355$/$01) is acknowledged.}\\
  \it Aston University, Birmingham B4 7ET, UK}
\begin{document}
\maketitle\vspace*{-7mm}
\baselineskip = 1.2\baselineskip
\begin{abstract}
Wavelet families arise by scaling and translations of a prototype 
function, called the {\em {mother wavelet}}. The construction 
of wavelet bases for cardinal spline spaces is generally  
carried out within the multi-resolution analysis scheme. 
Thus, the usual way of increasing the dimension of 
the multi-resolution subspaces is by augmenting  the scaling factor. 
We show here that, when working on a compact interval, 
the identical effect can be achieved without 
changing the wavelet scale but reducing the translation 
parameter. By such a procedure we generate a redundant frame, 
called a {\em{dictionary}}, spanning the same spaces as a 
wavelet basis but with wavelets of broader support. 
We characterise the correlation of the dictionary elements by 
measuring their `coherence' and produce examples illustrating 
the relevance of highly coherent dictionaries to problems of 
sparse signal representation. 
\end{abstract}
%\begin{keywords}
 % {\em Keywords:}
%  B-splines, wavelets, redundant dictionaries, non-linear approximation%in, sparse representation.
%\end{keywords}
\section{Introduction}
Canonical Coherent States have played an important 
role in quantum physics from the early days to the present time
\cite{Gla63,Dau80,KS85,Kai90,AAG91a,AAG00,Vou06,RAT06,DSA06}.
Mathematically they have been studied within the structure of
frames \cite{AAG93,Kai94,AAG00}. Their equivalent in 
the context of signal processing  
were first introduced by Gabor and  nowadays  
frequently appear in some contexts under the name of Gabor frames 
or Weyl-Heisenberg   
frames \cite{HW89,Dau92}.
On the other hand, the particular class of coherent states arising
from the affine group on the real line, called wavelets 
\cite{AK69,HW89,AAG00,Kai94,Dau92,Reb01}, 
have also been broadly applied in physics \cite{Kai03,
ABGK00,AK04,Van99}, as well as in signal processing.
The applications of wavelets notably increased  
when fast techniques,
arising from the multi-resolution analysis scheme, 
became available. 
Ever since wavelet analysis has been a popular tool 
among practitioners. 

For the most part multi-resolution analysis has been 
applied for generating orthogonal, semi-orthogonal and
biorthogonal wavelet bases.
Some redundant frames have been constructed  within 
this framework as well. However,
we should recall that 
countable sets of frame wavelets were first constructed  
by discretisation  of the translation and scaling parameters 
of affine coherent states.
%When such a discretisation process yields a redundant frame 
%in signal processing jargon 
%is referred to as 'oversampling'.
According to the mathematical 
measure of coherence introduced in \cite{Tro04}  
an orthogonal basis has coherence zero  and 
a non-orthogonal basis, although not redundant by definition, 
has coherence greater than zero. It is the purpose of the 
present effort to show that, 
in some finite dimensional spaces, by 
{\em {properly increasing the coherence   
 one can generate spaces of higher dimension}}. 

The success in constructing a number of linearly independent 
wavelets with good mathematical properties, such as
regularity or localisation, is very much appreciated by 
researchers in different fields (including the authors of 
this contribution) as a major achievement in
the design of wavelets. Nevertheless, in this communication we 
would like to illustrate, by recourse to a very interesting example,
a remarkable property of a class of highly coherent 
redundant wavelet systems. We will construct such systems for 
finite dimensional {\em{cardinal splines spaces}}, 
which are well characterised from a 
mathematical viewpoint \cite{Sch81} 
and will be shown to be useful for 
the problem of sparse signal representation by 
non-linear approximation. 

The problem of non-linear approximation concerns the
representation of a given function (signal) through 
the selection of waveforms, frequently called {\em{atoms}}, which 
are taken from a redundant set, called
a {\em{dictionary}} \cite{Mal98}. 
This problem has been the subject of
quite recent theoretical work with
regard to quasi incoherent dictionaries
\cite{GN04,Tro04}.
%Moreover, a number of techniques
%for selecting atoms have been devised in the
%signal processing field over the last fifteen years.
%\cite{chen:91,mallat:93,pati:93,davis:97,jaggi:98,chen:99,
%neira:02,andrle:04, andrle:04c}. 
%An advantage in the use of dictionaries for
%signal representation is often recognised
%by the possibility of considering atoms
%of different nature to match different structures of
%a signal. 
%However, the main result of the present effort is to
%show that, at least for cardinal spline spaces, one can
%construct useful dictionaries consisting of highly coherent 
%wavelets of the same nature as that of
%the wavelet basis for the identical space.
%This extends the property of cardinal B-spline dictionaries 
%demonstrated in \cite{andrle:04b}.

The term {\em {spline wavelets}} comprises a number of wavelet systems
ranging from Haar piece-wise constant wavelets to wavelet functions
of much higher regularity. In our framework all such systems admit an
equivalent construction and a particular one is obtained by
setting the order of the corresponding B-spline scaling function.
We will restrict wavelet systems to an interval and
construct spline wavelet dictionaries on the basis of the
following result: 
Let us focus on a cardinal spline space on a compact
interval with distance $2^{-j},\, j> 0$ between two adjacent knots.
We prove that such a space can be spanned by translating a wavelet
taken from a multi-resolution subspace corresponding to a
fixed scale, say scale $2^i$, $0<i<j$,
as long as the distance between two consecutive functions is
reduced to be $2^{i-j}$. This interesting feature provides the 
foundations for the
construction of a large variety of possible dictionaries for the
identical space.  In particular, 
multi-resolution-like dictionaries can be constructed
by spline wavelets whose support at the finest scale is
 larger than that corresponding to the finest scale in a
 multi-resolution analysis of the identical space. 
 This property provides a clear explanation of
 the `power of coherence'. Our results for finite dimensional 
 cardinal spline spaces
 on a compact interval assert that the
 benefit of increasing coherence, by decreasing the  
 translation parameter of the functions, is not only a consequence 
 of incorporating redundancy but also of the fact that by such a 
 procedure one may 
 {\em increase the dimension of the space}. This phenomenon 
  emerges clearly from our construction.

As will be illustrated
by numerical simulations, transforming a spline wavelet basis into
highly coherent dictionaries has a significant impact on 
signal representation by
non-linear techniques. It will be shown that some of  these 
dictionaries may yield a significant gain in the sparseness of a
representation with respect to the results produced by means of the
corresponding basis. This is enhanced by the comparison with other
techniques such as Best Basis Selection using Wavelet Packets.

%This outcome is very encouraging 
%because spline wavelet bases are known to be appropriate tools 
%for signal processing purposes \cite{unser:99}. Moreover, 
%the comparison 

The paper is organised as follows: Section~\ref{sec:multi}
introduces some background on spline multi-resolution analysis on a compact
interval relevant for our purpose. Section~\ref{sec:dic} establishes
the fact that, rather than increasing the wavelet scale to provide a
representation of a cardinal spline space of higher dimension, the
representation can be achieved by appropriate reduction of the
distance between two consecutive wavelets.  This result is used for
constructing a multi-resolution-like redundant dictionary that, through
the numerical examples of Section~\ref{sec:app}, is shown to be
useful for sparse signal representation.  The conclusions are drawn in
Section~\ref{sec:con}.

\section{Cardinal spline multi-resolution analysis on a compact
  interval}\label{sec:multi}
Let us recall that  multi-resolution analysis restricted to the
interval $[c,d]$, involves a sequence of nested spaces
% \begin{equation*}
$\V_0\subset \V_1\subset \cdots$
% \end{equation*}
satisfying that
% \begin{equation*}
$ \bigcup_{j\in\Z^+} \V_j\text{ is dense in } \L.$
% \end{equation*}
The complementary wavelet subspaces $\W_j$ are constructed in order to
fulfil
\begin{equation}
  \V_{j+1}=\V_j\oplus \W_j,\ j\in\Z^+
\end{equation}
so that
\begin{equation}
  \L=\V_0\oplus \W_0 \oplus \W_1\oplus \cdots = \V_0\oplus \bigoplus_{j\in\Z^+} \W_j.
\end{equation}
Without loss of generality we will assume throughout the paper that
$c,d \in \Z$.  Let us consider now that $\V_j, j\geq 0$ are cardinal
spline spaces of order $m$ with simple knots at the equidistant
partition of the interval $[c,d]$, having distance $2^{-j}$
between two adjacent knots. This implies that each $\V_j$ consists of piece-wise
polynomials of order $m$ having $m-2$ continuous derivatives
\cite{Sch81}.  A basis for $\V_j$ arises from the restriction of
the functions
\begin{align}\label{sca}
  \phi_{j,k}(x):=2^{j/2}\phi(2^j x-k), \ k\in\Z
\end{align}
to the interval $[c,d]$. The corresponding scaling function
$\phi(x)\equiv \phi_{0,0}(x)$ is the cardinal B-spline of order $m$
associated with the uniform simple knot sequence $0,1,\dots,m$.  Such
a function is given as
\begin{equation}
  \phi(x)=\frac{1}{m!}\sum_{i=0}^m(-1)^i\binom{m}{i}(x-i)^{m-1}_+,
\end{equation}
where $(x-i)^{m-1}_+$ is equal to $(x-i)^{m-1}$ if $x-i>0$ and 0
otherwise.

It should be stressed that different ways of constructing the boundary
functions give rise to different bases for $\V_j$. The 
restriction of the scaling functions  to the interval $[c,d]$
provides one  such basis, which is easy to construct.
Indeed, considering that the support of $\phi_{j,k}(x)$ is
$\supp\phi_{j,k}=\left[\frac{k}{2^j},\frac{k+m}{2^j}\right],$ one
has:
\begin{equation}
  \V_j=\Span\{\phi_{j,k} : k\in(2^jc-m,2^j d)\cap \Z\}. \label{basis}
\end{equation}
%Since the dimension of a spline space is the number of inner knots
%plus the order of the splines we have: $\dim{\V_j}=(d-c)2^j-1+m$.  As
%required, this is the cardinality of the index set in \eqref{basis}.
Without loss of generality, we assume that at least one scaling
function from $\V_0$ and one wavelet from $\W_0$ are completely
contained in $[c,d]$. This is equivalent to assuming that $d-c\geq
\max\{m,w\}$ where $w$ is the length of $\supp \psi$.  A basis for
$\W_j$ can then be constructed by restricting the functions
\begin{equation}\label{psi}
 \psi_{j,k}(x):=2^{j/2}\psi(2^j x-k),\,  k\in \Z
\end{equation}
to the interval $[c,d]$ and eliminating some redundancy introduced by
the cutting process.  One can start by considering all the functions
$\psi_{j,k}(x)$ having non-trivial intersection with the interval
$(c,d)$, which restricts the values of the index $k$ to
$k\in(2^jc-w,2^j d)\cap \Z$.  Since the cardinality of such a set is
$(d-c)2^j-1+w$, but $\dim \W_j=\dim \V_{j+1}-\dim \V_{j}=(d-c)2^j$, in
order to have a basis it is necessary to eliminate $w-1$ boundary
functions.  A natural choice is to eliminate the first $\lceil
z\rceil$ left boundary functions and the last $\lfloor z \rfloor$
right boundary functions where $z=(w-1)/2$ and $\lceil\cdot\rceil$,
$\lfloor \cdot\rfloor$ indicate the upper and lower integer part,
respectively.  A basis for $\W_j$ constructed by the cut-off process
is thereby given as:
\begin{equation}
  \W_j=\Span\left\{\psi_{j,k} :
    k\in \left(2^jc- \left\lfloor z\right\rfloor,
      2^jd- \left\lfloor z\right\rfloor \right)\cap \Z\right\}.
  \label{basis2}
\end{equation}
\section{Spline wavelet dictionaries}\label{sec:dic}
In a previous publication we have shown that a cardinal spline space on
a compact interval can be spanned by dictionaries consisting of
functions of broader support than the corresponding B-spline basis
\cite{AR05}. In the present notation this entails that for an
integer $\ell\geq 1$ the set of B-splines of order $m$
\begin{equation}\label{sdic}
  \dicv_{j,\ell}=\{\phi_{j,\kj}:
  \kj\in(2^jc-m,2^jd)\cap\Z/2^{\ell} \}
\end{equation}
satisfies
 \begin{equation}\label{teo0}
 \Span\{\dicv_{j,\ell}\}=\V_{{j+\ell}}.
 \end{equation}
Note that the support of a function in \eqref{sdic} completely
contained in $[c,d]$ is $2^\ell$-times  the support of a
corresponding basis function from $\V_{{j+\ell}}$, and we have
introduced a new translation parameter $\kj$ which is no longer an
integer but $\kj\in\Z/2^{\ell}$.

We call the set $\dicv_{j,\ell}$ a {\em{B-spline dictionary}} and wish
to consider now the possibility of creating cardinal spline
dictionaries of wavelets. To this end we construct the set
$\dicw_{j,\ell}, \, \ell\geq 1$ as follows: The wavelets
are translated using the translation parameter, $\kj$,
defined above. All the functions $\psi_{j,\kj}(x), \kj\in\Z/2^{\ell}$
whose supports have non-trivial intersection with the interval $(c,d)$
are considered.  Thus
\begin{equation}{\label{wdic}}
  \dicw_{j,\ell}:=\{\psi_{j,\kj}
  : \kj\in (2^jc-w,2^jd)\cap \Z/2^{\ell}\}.
\end{equation}
Our proposal for constructing cardinal spline wavelet dictionaries
stems from the fact that $\Span\{\dicw_{j,\ell}\}=\V_{{j+\ell}}$. This
result is established in the next theorem. The proof is achieved  by 
 using an equivalent technique to that for the proof of Theorem~1  
in \cite{AR05} and is given in
Appendix A.
\begin{theorem}\label{main_theo}
  Let $\dicw_{j,\ell}$ be given as in \eqref{wdic} for a given integer
  $\ell\geq 1$.  Then the following relation holds:
  \begin{equation}
    \Span\{\dicw_{j,\ell}\}=\V_{{j+\ell}}.
  \end{equation}
\end{theorem} 
\begin{remark}\label{re1}
The above theorem gives us a tool for designing different wavelet
dictionaries for cardinal spline spaces. Notice that
simply by setting different values of the index $\ell$ in
\eqref{wdic} we obtain dictionaries of wavelets of different support
spanning the identical cardinal spline space $\V_{{j+\ell}}$. 
\end{remark}

As described below the result of Theorem \ref{main_theo}
allows us to create multi-resolution-like
dictionaries for the space being considered.

\begin{figure}[!t]
  \begin{center}
   \includegraphics[width=8cm]{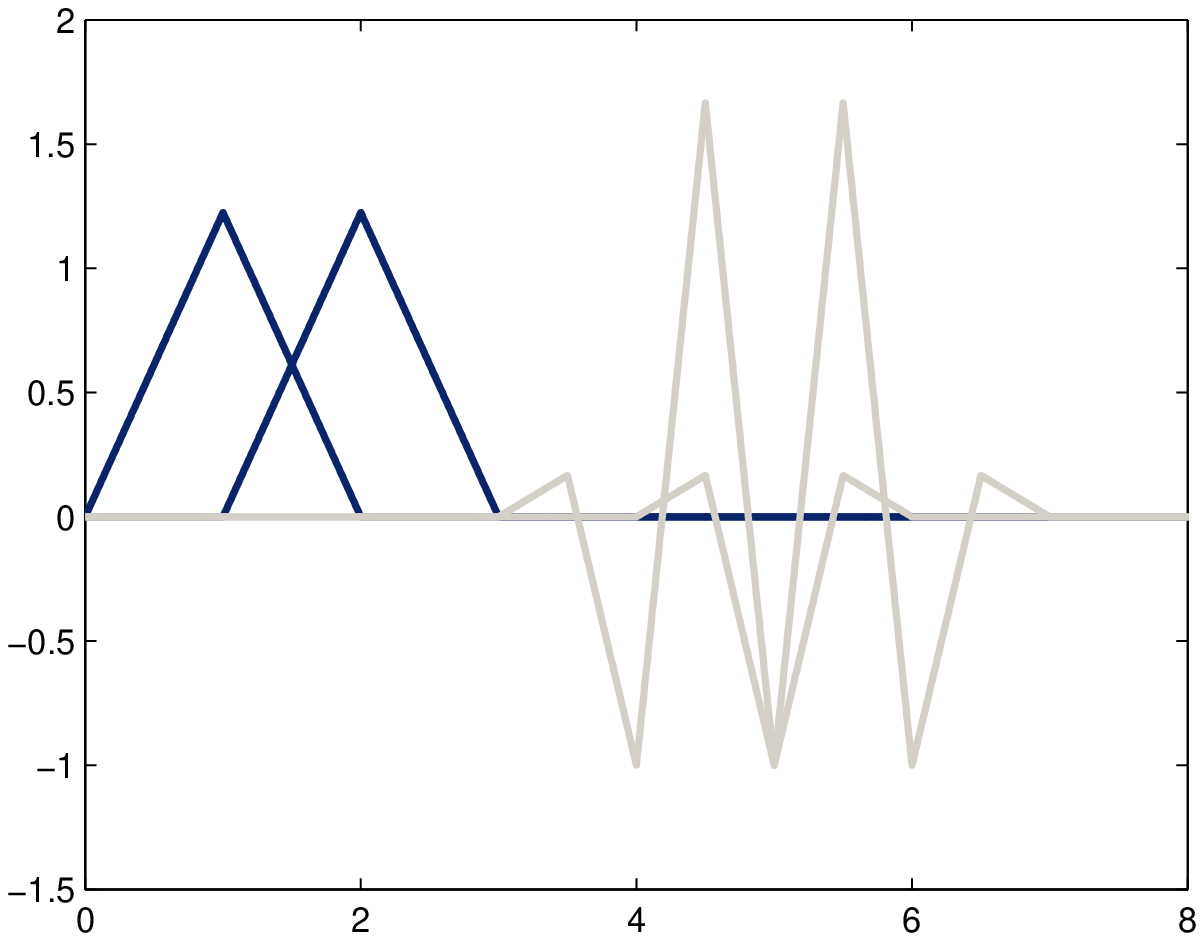}%\\ \vspace*{-3mm}
   \includegraphics[width=8cm]{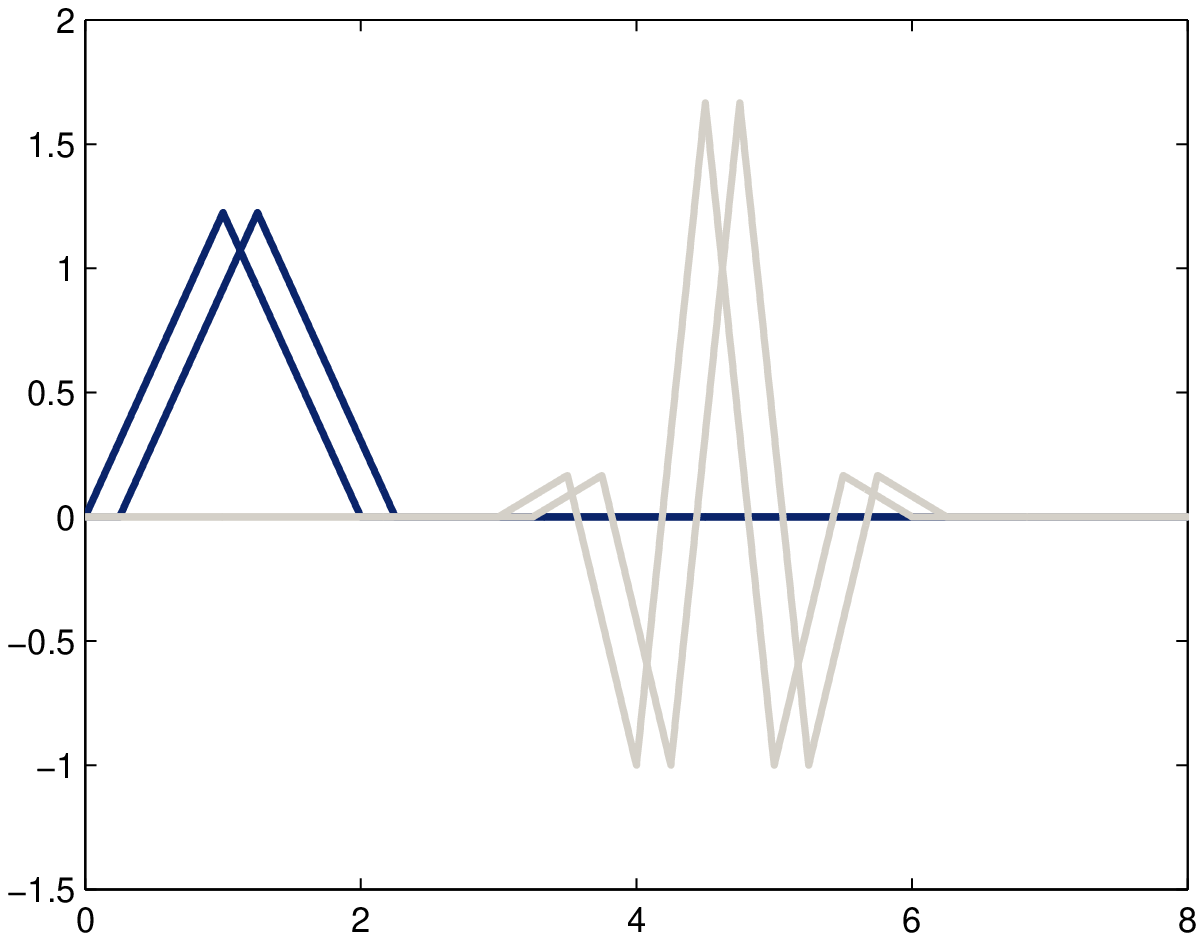}\\
   \includegraphics[width=8cm]{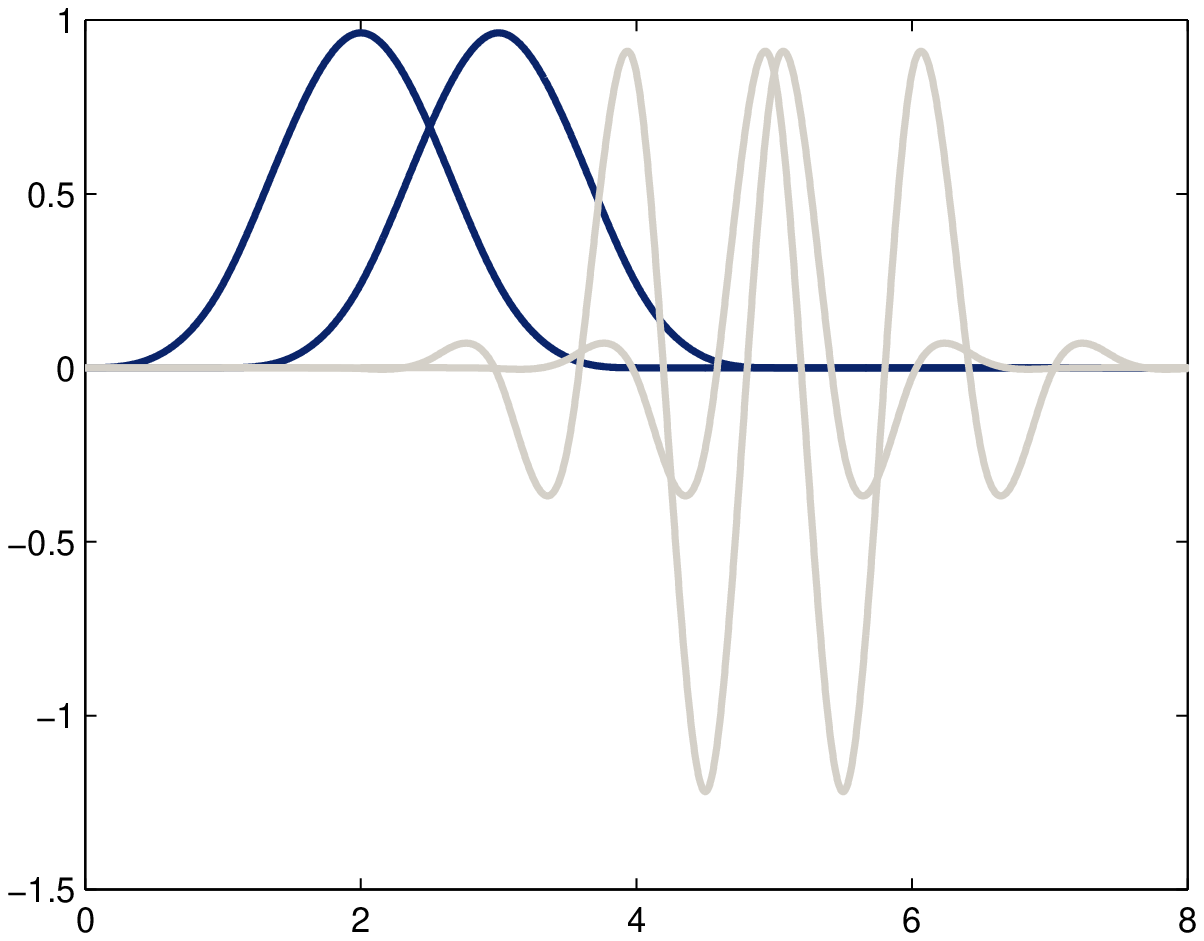}%\\ \vspace*{-3mm}
   \includegraphics[width=8cm]{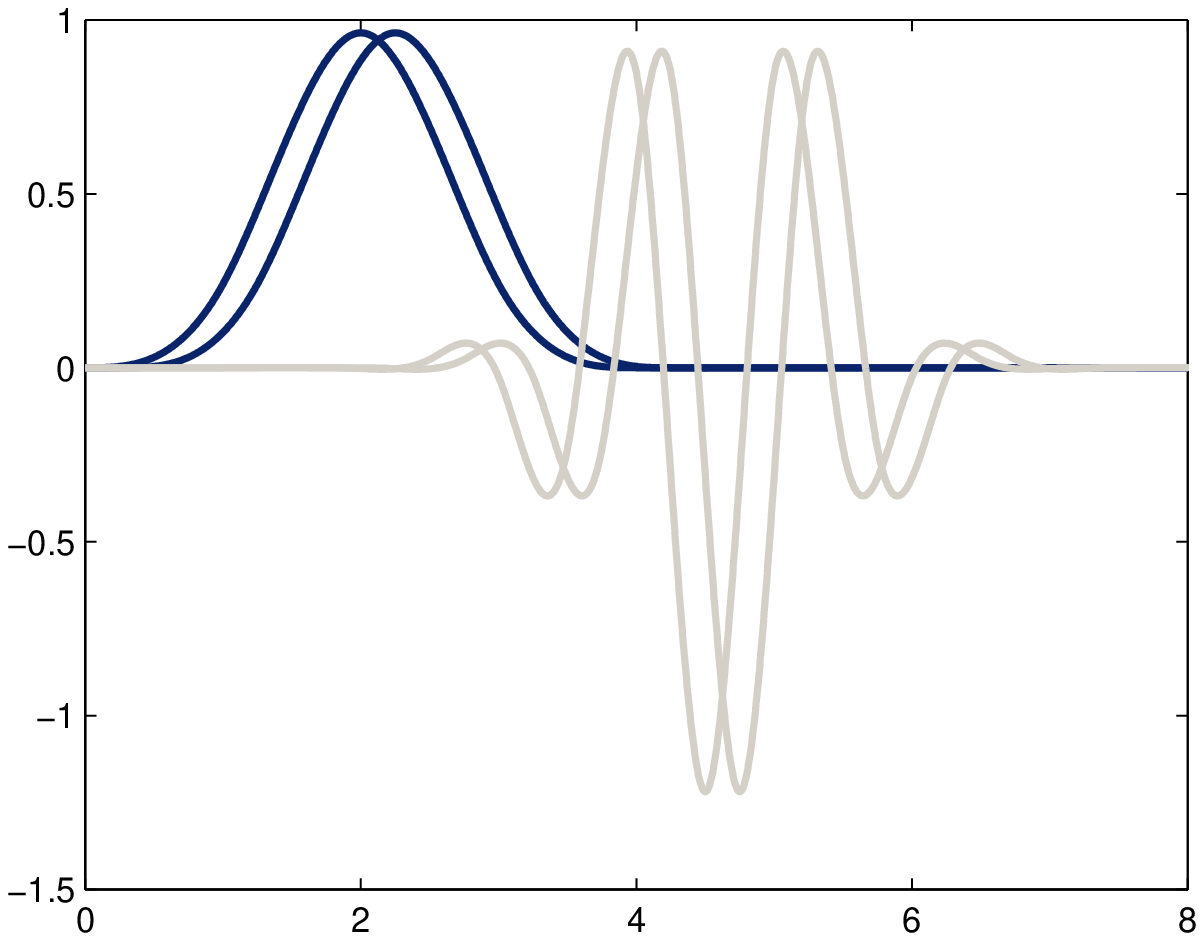}

  \end{center}\vspace{-11mm}
  \caption{The top left graph shows two
  consecutive scaling functions (dark line) and  two
  consecutive wavelets at the coarsest scale (light line) in a 
  linear  spline wavelet basis (order $m=2$). 
  The right graph shows the equivalent functions in a dictionary
  constructed by considering $\ell=2$.
  The bottom graphs have the same description as the top ones, 
  but involving
  cubic spline wavelets (order $m=4$).}
 \label{fd1}\vspace*{7mm}
\end{figure}

We denote by $\dicv_{j,0}$ the B-spline basis \eqref{basis} for
$\V_{j}$ and by $\dicw_{j,0}$ the spline wavelet basis \eqref{basis2}
for $\W_{j}$.  The classical wavelet decomposition $\V_{j}=
\V_0\oplus\W_0\oplus\cdots\oplus\W_{j-1}$ in terms of bases can then
be expressed in our ``dictionary notation'' as $\V_{j}=\Span
\{\D_{j,0}\}$ where:
\begin{equation}\label{bas}
  \D_{j,0}=\dicv_{0,0}\cup \dicw_{0,0} \cup \dicw_{1,0} \cup \cdots 
  \cup \dicw_{j-1,0}.
\end{equation}
Considering $\ell\geq 1$ we construct multi-resolution-like
dictionaries spanning $\V_{{j}}$ as
\begin{equation}
  \label{dic}
  \D_{j,\ell}=\dicv_{0,\ell}\cup \dicw_{0,\ell} \cup \dicw_{1,\ell} \cup \cdots 
  \cup \dicw_{j-\ell,\ell}.
  % \D_{{j},\ell}=\dicv_{0,\ell}\cup\bigcup_{i=0}^{j-\ell}
  % \dicw_{i,\ell}.
\end{equation}  
\begin{remark}
It follows from \eqref{bas} and \eqref{dic} that the minimum support 
of an inner wavelet in $\D_{{j},\ell}$ is
$2^{\ell-1}$-times the minimum support of an inner wavelet
in $\D_{j,0}$.
\end{remark}
%This is advantageous because it implies that 
%the dictionary functions could
%be numerically represented by using a coarser discretisation of the
%continuous variable.
The top left graph of Fig.~\ref{fd1} shows two
consecutive scaling functions (dark line) and  two 
consecutive wavelets at the coarsest scale (light line) in a 
linear spline basis  (order $m=2$).  
The right graph shows the equivalent functions in a dictionary 
constructed by considering $\ell=2$.
The bottom graphs depict an equivalent example but involving 
cubic splines, which correspond to order $m=4$. Let us remark that,  
since in  both examples the dictionaries arise by 
setting $\ell=2$, for representing the same space as 
the corresponding basis we  should  consider a coarser scale, i.e, 
if the finest scale in the basis corresponds to a scaling 
factor $2^{j}$ in  the
 dictionary the finest scale corresponds to a scaling factor 
 $2^{j-1}$. 

Of course, considering the result of Theorem~\ref{main_theo} 
one could construct many different dictionaries for the 
identical space. 
Each such dictionary constitutes a frame for a cardinal spline 
space of order $m$ with distance $2^{-j}$ between knots.  As
will be illustrated in the next section, particular dictionaries
constructed as in \eqref{dic} may yield a significant gain, as far as
sparseness is concerned, in problems of signal representation. 

For the examples given in Section~\ref{sec:app} we construct specific
dictionaries on the interval $[0,8]$ as follows:
We consider cubic splines
and use the semi-orthogonal wavelets (Chui-Wang4 family) introduced in
\cite{CW92} to construct the basis for each $\dicw_{i,0},\,
i=1,\ldots,j-1$ by the
simple cut off process described in Section~\ref{sec:multi}.
In order to construct the dictionary 
we take a prototype function 
from each subspace $\dicw_{i,0},\,i=1,\ldots,j-2$ 
(and also from $\dicv_{0,0}$)
and translate such a function  to a distance $2^{-2}$. Notice that 
we do not use functions from $\dicw_{j-1,0}$. This is because 
by reducing the distance between functions to $2^{-2}$ we 
need one less scale than the wavelet basis to represent the 
same space.  
%In doing so, due to the
%boundary functions, the semi-orthogonality property (orthogonality 
%between different scales) of the
%decomposition is lost. Since the support of $\psi(x)$ is of 
%length $w=2m-1=7$ if one wishes to restore semi-orthogonality
%one must orthogonalize  $\lceil
%(w-1)/2\rceil=3$ left boundary functions and  $\lfloor (w-1)/2
%\rfloor=3$ right boundary functions with respect to $V_0$ and 
%so on.
%For the only sake of comparing coherence we will consider both 
%constructions on $[0,8]$ and
%denote as $\Dsech_{j}$ the semi-orthogonal decomposition 
%of $V_{j}$. 

\begin{figure}[!t]
\begin{center}
 \includegraphics[width=6cm]{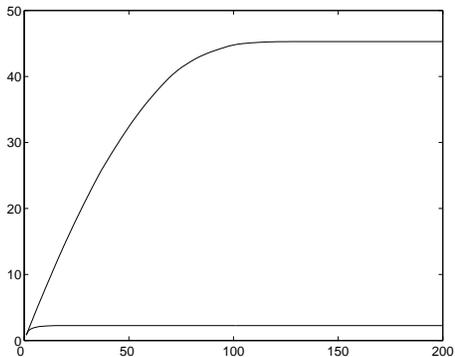}
%   \end{center}\vspace{-11mm}
\caption{Cumulative coherence function $\mu(p),\,p=1,\ldots,200$
as defined in \eqref{cuco}.
The lower line corresponds to the basis $\D_{6,0}$
of cubic spline wavelets for $V_6$. The other line corresponds to the 
dictionary $\D_{6,2}$ for the same space.}
\label{coe}
\end{center}
\vspace{9mm}
\end{figure}

In order to characterise the coherence of the dictionaries at hand
we use the 
{\em {cumulative coherence}} function, $\mu(p)$, introduced 
in \cite{Tro04}. Given a dictionary 
$\{\alpha_{\omega}\}_{\omega \in \Omega}$, where
 $\Omega$ is the set of indices labelling the dictionary atoms, 
 $\mu(p)$  measures how much a collection of 
$p$ atoms resembles a fixed one. It is defined by
\begin{equation} \label{cuco}
\mu(p)=  \max_{|\Lambda| =p} \max_{\omega \not\in \Lambda} 
\sum_{\lambda \in \Lambda}|\la \alpha_\omega,  \alpha_\lambda\ra|,
\end{equation}
where $|\Lambda|$ indicates the cardinality of the set $\Lambda$.

The cumulative coherence
for the basis  $\Dch_{6,0}$ and the dictionary $\Dch_{6,2}$, 
which will be used in the next section, 
is plotted in  Fig.~\ref{coe}. It is clear from this figure 
that the coherence of the dictionary $\Dch_{6,2}$  
is much more larger than the coherence  of the basis 
$\Dch_{6,0}$. 
\begin{remark}
An increment in coherence also implies that 
the construction of the  dual  
frame for the dictionary may 
become an ill posed problem. However, when the aim is to use 
dictionaries for sparse signal representation 
one is not interested in the dual frame for the whole 
space. On the contrary, in order to 
produce orthogonal projections onto the approximation subspace, 
the duals for the particular  subspace  need to be determined.  
For the representation to be useful for compression purposes the 
dimension of the approximation subspace should be considerable 
smaller that the dimension of the  dictionary space. 
Thus,  
the construction of the duals in the corresponding subspace
is well posed. 
Recursive techniques for 
 adapting dual functions  are discussed in 
\cite{Reb02,Reb04b,Reb07a}.
\end{remark}

\section{Application to sparse signal representation}\label{sec:app}
\begin{figure}[!t]
  \begin{center}
    \includegraphics[width=8cm]{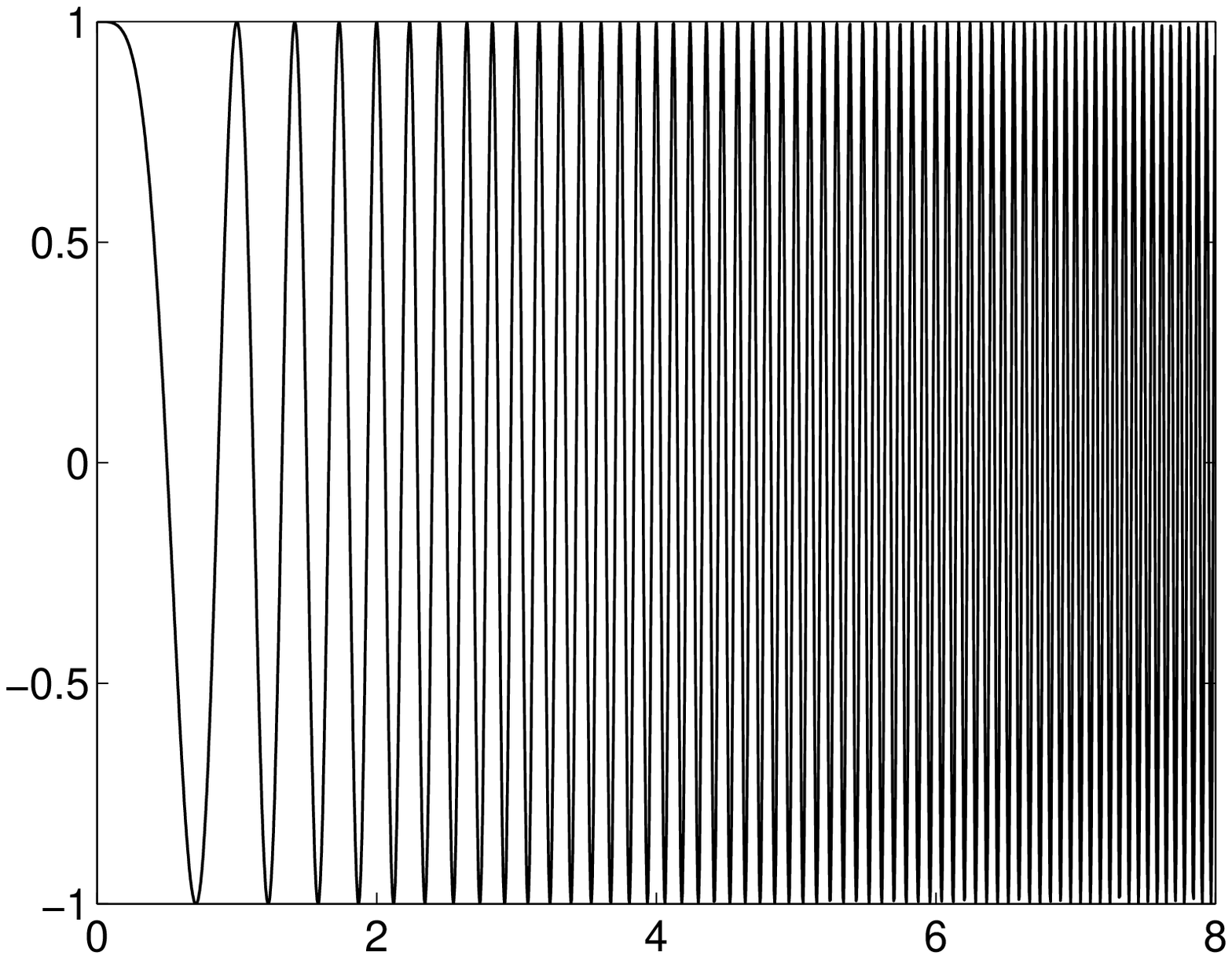}%\\ \vspace*{-3mm}
    \includegraphics[width=8cm]{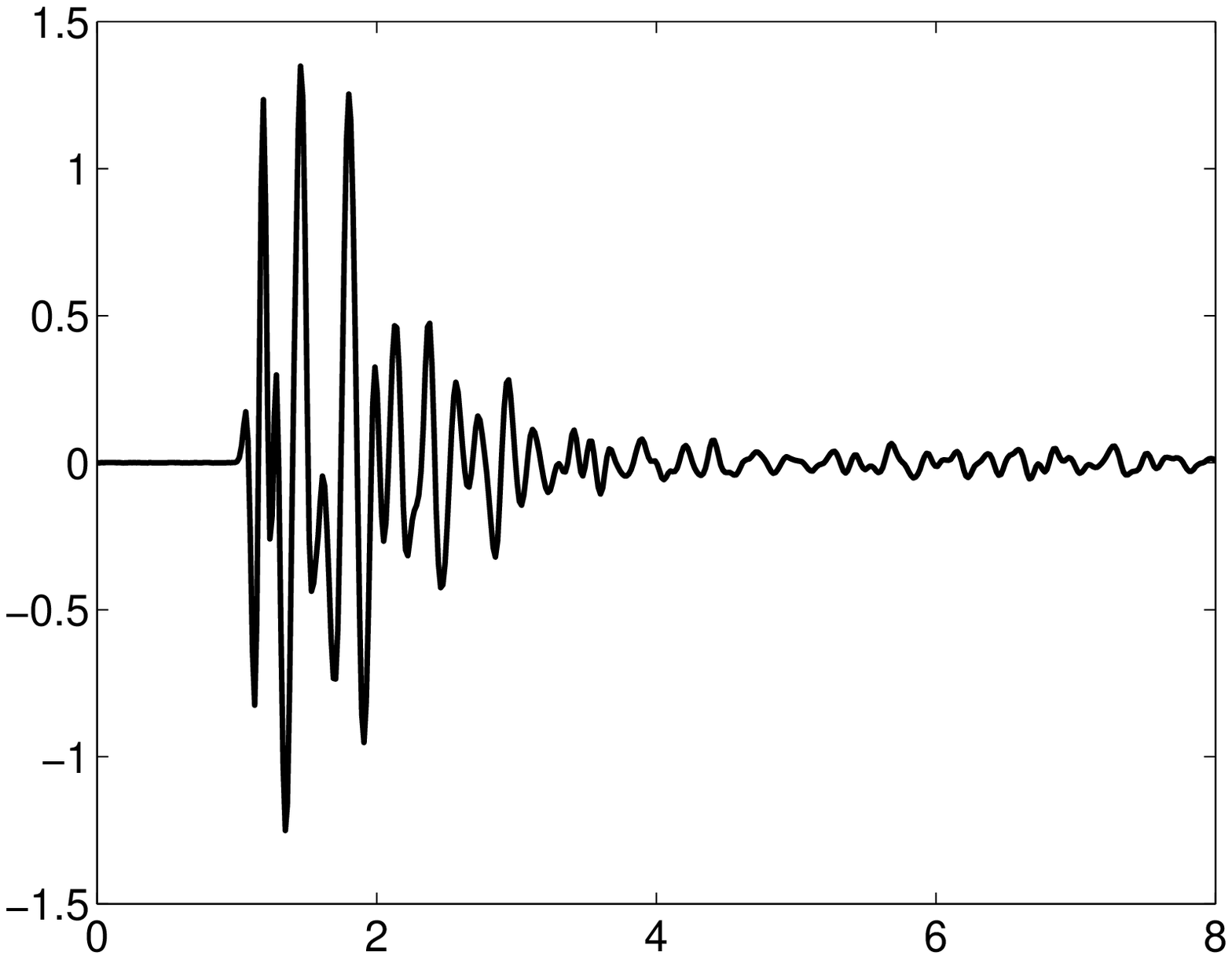}
  \end{center}
 % \vspace{-11mm}
  \caption{Chirp signal $f=\cos(2\pi x^2)$ (left). Seismic
    signal (right).}
  \label{fig1}
  \vspace{7mm}
\end{figure}

In this section we use the multi-resolution-like
dictionary characterised in the previous section, 
corresponding to considering $\ell=2$. 
The central aim is to compare the sparseness of the 
signal representation
achieved by using the basis $\Dch_{j,0}$ and 
the dictionary $\Dch_{j,2}$.
Moreover, the results produced by a number of other
techniques will be presented for further comparison.
%
% such as Best Basis (BBS) selection from {\em Wavelet Packets}
% \cite{coifman:92,wick:94} will be also presented for further
% comparison.
%
The two signals to be represented are the chirp and the seismic signal of Fig.~\ref{fig1}.

The  piece of seismic signal
was taken from the WaveLab802 Toolbox \cite{wavelab} (it is
acknowledged there that such a signal is distributed throughout the
seismic industry as a test dataset). 
Both signals have a good approximation, coinciding with the graphs
of Fig.~\ref{fig1}, in $V_6$, the space of cardinal cubic splines 
with the distance $2^{-6}$ between two adjacent knots. 
In order to represent this space we will use:
The basis $\Dch_{6,0}$ and the dictionary
$\Dch_{6,2}$. 
%It is important to remark that the wavelets of the 
%finest scale in $\Dch_{6,0}$ are {\em not present} in $\Dch_{6,2}$, as
%they are not needed when the required space $V_6$ is spanned by this
%dictionary.  
In both cases we apply the same strategy for selecting
the atoms $\{\alpha_{l_n}\}_{n=1}^N$ to
approximate the given signal by the {\em atomic decomposition}
\begin{equation}\label{atode}
 f^N= \sum_{n=1}^N c_n^N \alpha_{l_n},
\end{equation}
where the superscript indicates that the coefficients $c_n^N$ yield
the orthogonal projection of the given signal onto the subspace
spanned by the selected atoms $\{\alpha_{l_n}\}_{n=1}^N$.

Our strategy evolves by stepwise minimum residual error selection and
is carried out in three stages:
\begin{itemize}
\item[i)]
The atoms are selected one by one according
to the adaptive pursuit method
\cite{RL02} in order to reach a predetermined precision in the
representation of the signal. 
\item[ii)]The previous approximation is improved by
means of a ``swapping procedure'' which operates as follows: at each
step an atom of the atomic decomposition is replaced by a dictionary
atom, provided that the operation improves the residual error
\cite{AR06}. 
\item[iii)]
A backward pursuit method \cite{ARS04} is applied to 
disregard some
coefficients  in order to produce an approximation up to the error 
of ~i).
\end{itemize}
MATLAB codes for constructing the proposed dictionaries and
implementing the selection strategy are available at \cite{AR03}.

\begin{table}[!t]
\caption{Comparison of the sparseness achieved by
 selecting atoms from the basis $\Dch_{6,0}$, the dictionary 
 $\Dch_{6,2}$, FWT and BBS from 
 WP using different criteria and wavelet 
 families.} 
\label{tab1}
\centering\vspace*{5mm}
 \begin{tabular}{|c|c|c|c|}\hline
    Family &Approach&  $N$ (chirp)& $N$ (seismic)\\
    \hline \hline
    Chui-Wang4 & Basis $\Dch_{6,0}$  & 254       &210      \\
    & Dictionary $\Dch_{6,2}$ & \bf 166   &\bf 125\\ 
    \hline
    Bior4.4  & FWT &  423 &223\\
    & BBS Shannon entropy&  315       &217\\
    & BBS $L^1$ norm   &  301       &206\\
    & BBS $L^{0.5}$ norm&  289       &206\\
    \hline
    Daub10  & FWT  &  401     &222\\
    & BBS Shannon entropy&  279       &201\\
    & BBS $L^1$ norm   &  271       &195\\
    & BBS $L^{0.5}$ norm&  271       &194\\
    \hline
  \end{tabular}
%  \vspace*{-3mm}
\end{table} 

For further comparison we have represented the signals using Fast
Wavelet Transform (FWT) with the biorthogonal (cubic spline) family
Bior4.4 and the Daubechies orthogonal family Daub10.  With these two
families we also considered a number of different criteria for Best
Basis Selection (BBS) from the corresponding Wavelet Packets (WP)
\cite{Coi92}.
These techniques were applied by using the Wavekit tool
\cite{wavekit} and iteratively tuning the threshold in order to
produce a signal approximation of the same quality (with respect of
the $L^2$ error norm) by all the approaches. 
The results are presented in
Table~\ref{tab1}. The first column labels the wavelet families. 
The second one specifies the approach and the third and fourth
column the number of atoms  needed by those approaches to represent 
the corresponding signal up to the same precision.

% We should point that our selection techniques are independent of the
% atoms norm and therefore normalisation is not necessary.
As  can be seen in Table~\ref{tab1}, to obtain the desired
approximation of the chirp of Fig.~\ref{fig1} using the basis
$\Dch_{6,0}$ we need 254 atoms.  Sparseness is improved to 166 atoms
by using the proposed redundant dictionary $\Dch_{6,2}$ involving
functions of the same nature.  On the contrary, the number of atoms
that are needed to represent the chirp signal by all the other
approaches is larger.  The same feature is exhibited in the
representation of the seismic signal. 
However, there is an interesting
phenomenon that can be observed: In the previous case the FWT 
(for both the
Bior4.4 and Daub10 families) was shown to produce a much poorer
sparseness property than the basis $\Dch_{6,0}$. In this case the
results yielded by the three wavelet bases are comparable. This is due
to the fact that the nature of the chirp near the left boundary is
very different from that near the right one.  Hence, the periodic
boundary conditions in the implementation of FWT have a very negative
effect in relation to the spareness of the representation. We believe
this is the reason why the cut-off approach performs much better in
representing the chirp.  However, in the case of the seismic signal
for which the two ends of the signal are not so different and of small
magnitude, the approximations obtained by all the wavelet bases are
comparable.  The representation of this signal by selecting atoms from
the proposed dictionary is again significantly superior to all the
other techniques.

\section{Conclusions}
\label{sec:con}
The construction of wavelet dictionaries for cardinal spline spaces on
a compact interval was discussed. It was first proved that if 
$2^{-j}$ is the distance between consecutive knots in a 
cardinal spline space on a compact interval, such  a
space can be spanned by translating wavelets at scale
$2^i$ a distance $2^{i-j}$. 
This result can be used for
building a large variety of wavelet dictionaries for the identical
space. In particular, multi-resolution-like  highly coherent
dictionaries of
wavelets of different supports were constructed. An interesting
feature of the proposed construction is that it allows us to span a
spline space of given dimension by using wider wavelets than those
necessary in multi-resolution analysis. As illustrated by the
numerical examples, this feature is relevant to sparse signal
representation by non-linear techniques. The simulations presented
here clearly show that the proposed multiresolution-like dictionaries,
which are easily constructed by simple translations of prototype
functions at different scales, may yield a very significant gain in
the sparseness of a signal representation by step-wise selection
techniques. This outcome could be somewhat `expected', since by 
 decreasing the translation step one also generates 
more choice so as to chose the suitable functions for representing 
a given signal. 
However, what is definitely a very interesting and somewhat  
surprising result is that, 
when working on a compact interval, 
by decreasing the translation step of a  wavelet family one can 
generate higher dimensional spaces. In our mind this is a
mathematical characterisation of the `power of coherence'. 

\subsection*{Acknowledgements}
Support from EPSRC (grants GR$/$R86355$/$01 and 
EP$/$D062632$/$1) is acknowledged.

\appendix
\section*{Appendix A: Proof of Theorem 1}
\begin{proof}
  The inclusion $\Span\{\dicw_{j,\ell}\}\subset\V_{{j+\ell}}$ is
  obvious because for $\ell\geq 1$ all the spline wavelets in $\W_j$
  are, by definition, in $\V_{{j+\ell}}$. Since the distance between
  two adjacent knots in $\V_{{j+\ell}}$ is $2^{-(j+\ell)}$, by moving
  the wavelets such a distance it is ensured that they remain in
  $\V_{{j+\ell}}$.  We can then build scaling-like equations for all
  the functions of the dictionary $\dicw_{j,\ell}$:
\begin{subequations}
    \begin{align}
      \psi_{j,\kj}(x)&=\hspace{-6mm}\sum_{n=2^{{j+\ell}}c-m+1}^{2^{\ell}(\kj+w)-m}\hspace{-6mm}
      g_{\kj,n}\phi_{{j+\ell},n}(x),\;   \kj\in (2^jc-w,2^jc)\cap
      \Z/2^{\ell},
      \label{1st}\\
      \psi_{j,\kj}(x)&=\hspace{-5mm}\sum_{n=2^{\ell}\kj}^{2^{\ell}(\kj+w)-m}\hspace{-5mm}
      g_{\kj,n}\phi_{{j+\ell},n}(x),\;  \kj\in [2^jc, 2^jd-w]\cap
      \Z/2^{\ell},
      \label{2nd}\\
    \psi_{j,\kj}(x)&=\hspace{-2mm}\sum_{n=2^{\ell}\kj}^{2^{{j+\ell}}d-1}\hspace{-2mm}
      g_{\kj,n}\phi_{{j+\ell},n}(x),\;  \kj\in (2^jd-w,2^jd)\cap
      \Z/2^{\ell}.
      \label{3rd}
    \end{align}
  \end{subequations}
  The proof of the inclusion
  $\Span\{\dicw_{j,\ell}\}\supset\V_{{j+\ell}}$ is achieved by showing
  that every function $\phi_{{j+\ell},n}(x) \in \V_{{j+\ell}}$ can be
  written as a linear combination of functions from $\dicw_{j,\ell}$.
  For every $ n\in[2^{{j+\ell}}c,2^{{j+\ell}}d)\cap \Z$ (corresponding
  to an inner or right boundary function $\phi_{{j+\ell},n}(x)$) there is a unique index
  $\kl=n/2^{\ell}$ such that $g_{\kr,n}\ne 0$ and $g_{\kr,i}=0$ for
  every $i<n$. Hence, in the line of \cite{AR05}, by using
  \eqref{2nd} or \eqref{3rd} the function $\phi_{{j+\ell},n}(x)$ is
  substituted by $\psi_{j,\kr}(x)/g_{\kr,n}$ plus a linear combination
  of $\phi_{{j+\ell},l}(x)$, $l>n$. This recursive process ends with
  the substitution $\phi_{{j+\ell},n}(x)=\psi_{j,\kr}(x)/g_{\kr,n}$
  where $n=2^{{j+\ell}}d-1$.

  As for every $ n\in(2^{{j+\ell}}c-m,2^{{j+\ell}}c)\cap \Z$
  (corresponding to a left boundary function $\phi_{{j+\ell},n}(x)$) there is also a
  unique index $\kl=(n+m)/2^{\ell}-w$ such that $g_{\kl,n}\ne 0$, and
  $g_{\kl,i}=0$ for every $i>n$.  Thus using \eqref{1st} the function
  $\phi_{{j+\ell},n}(x)$ is substituted by $\psi_{j,\kl}(x)/g_{\kl,n}$
  plus a linear combination of $\phi_{{j+\ell},l}(x), l<n$.  The
  recursive process ends with the substitution
  $\phi_{{j+\ell},n}(x)=\psi_{j,\kl}(x)/g_{\kl,n}$ where
  $n=2^{{j+\ell}}c-m+1$.

  Consequently, by back-substitution, we have the decomposition of
  each function $\phi_{{j+\ell},n}(x),\,
  n\in(2^{{j+\ell}}c-m,2^{{j+\ell}}d)\cap \Z$ in terms of
  $\psi_{j,\kj}(x),\ \kj\in (2^jc-w,2^jd)\cap\Z/2^{\ell}$.
\end{proof}

%\bibliographystyle{IEEEtran} 
%\bibliography{biblio,biblio-dict}
\bibliographystyle{elsart-num}
\bibliography{signal2}

\end{document}